\documentclass[reqno,centertags, 12pt]{amsart}
\usepackage{amsmath,amsthm,amscd,amssymb}
\usepackage{latexsym}

\newcommand{\bbN}{{\mathbb{N}}}
\newcommand{\bbR}{{\mathbb{R}}}

\newcommand{\bbC}{{\mathbb{C}}}

\newcommand{\calC}{{\mathcal C}}
\newcommand{\calL}{{\mathcal L}}
\newcommand{\calM}{{\mathcal M}}

\newcommand{\dott}{\,\cdot\,}
\newcommand{\no}{\nonumber}
\newcommand{\lb}{\label}
\newcommand{\f}{\frac}

\newcommand{\ol}{\overline}

\newcommand{\loc}{\text{\rm{loc}}}

\newcommand{\dom}{\text{\rm{dom}}}

\newcommand{\supp}{\text{\rm{supp}}}
\newcommand{\AC}{\text{\rm{AC}}}
\newcommand{\sabc}{\text{\rm{s.-a.\,b.c.}}}
\newcommand{\bi}{\bibitem}

\newcommand{\beq}{\begin{equation}}
\newcommand{\eeq}{\end{equation}}
\newcommand{\ba}{\begin{align}}
\newcommand{\ea}{\end{align}}
\newcommand{\veps}{\varepsilon}

\DeclareMathOperator{\Real}{Re}
\DeclareMathOperator{\Ima}{Im}

\DeclareMathOperator*{\eqlim}{=}

\allowdisplaybreaks
\numberwithin{equation}{section}

\newtheorem{theorem}{Theorem}[section]

\newtheorem{lemma}[theorem]{Lemma}
\newtheorem{corollary}[theorem]{Corollary}
\theoremstyle{definition}

\theoremstyle{remark}
\newtheorem{remark}[theorem]{Remark}

\begin{document}
\title[Local Borg-Marchenko Results]{On Local Borg-Marchenko 
Uniqueness Results}
\author[F. Gesztesy and B. Simon]{Fritz Gesztesy$^{1}$ and 
Barry Simon$^{2}$}
\date{October 15, 1999}

\footnotetext[1]{Department of Mathematics, University of 
Missouri, 
Columbia, MO 65211, USA. E-mail: fritz@math.missouri.edu}
\footnotetext[2]{Division of Physics, Mathematics, and 
Astronomy, 
253-37, California Institute of Technology, Pasadena, 
CA~91125, USA. 
E-mail: bsimon@caltech.edu. This material is based upon work 
supported by the National Science Foundation under Grant 
No.~DMS-9707661. The Government has certain rights in this 
material.}

\begin{abstract} We provide a new short proof of the 
following 
fact, first 
proved by one of us in 1998: If two Weyl-Titchmarsh 
$m$-functions, 
$m_j(z)$, 
of two Schr\"odinger operators $H_j = -\f{d^2}{dx^2} + q_j$, 
$j=1,2$ in $L^2 
((0,R))$, $0<R\leq \infty$, are exponentially close, that is, 
$|m_1(z)- m_2(z)| \underset{|z|\to\infty}{=} 
O(e^{-2\Ima (z^{1/2})a})$,
$0<a<R$, then
$q_1 = q_2$  a.e.~on $[0,a]$. The result applies to any 
boundary
conditions at
$x=0$ and 
$x=R$ and should be considered a local version of the 
celebrated 
Borg-Marchenko 
uniqueness result (which is quickly recovered as a 
corollary to 
our proof). 
Moreover, we extend the local uniqueness result to 
matrix-valued 
Schr\"odinger 
operators. 
\end{abstract}

\maketitle

\section{Introduction} \lb{s1}

Let $H_j = -\f{d^2}{dx^2} + q_j$, $q_j\in L^1 ([0,R])$ for 
all $R>0$, $q_j$ 
real-valued, $j=1,2$, be two self-adjoint operators in 
$L^2 ([0,\infty))$ with a 
Dirichlet boundary condition at $x=0_+$. Let $m_j(z)$, 
$z\in\bbC\backslash\bbR$ be 
the Weyl-Titchmarsh $m$-functions associated with $H_j$, 
$j=1,2$. The principal 
purpose of this note is to provide a short proof of the 
following 
uniqueness 
theorem in the spectral theory of one-dimensional Schr\"odinger 
operators, 
originally obtained by Simon \cite{Si98} in 1998. (Actually, 
Simon's result 
\cite{Si98} was weaker; the result as stated is from 
\cite{GS98}.)

\begin{theorem}\lb{t1.1} Let $a>0$, $0<\veps<\pi/2$ and 
suppose that
\begin{equation} \lb{1.1}
|m_1 (z) - m_2(z)| \eqlim_{|z|\to\infty} 
O(e^{-2\Ima (z^{1/2})a}) 
\end{equation} 
along the ray $\arg(z) = \pi-\veps$. Then
\begin{equation} \lb{1.2}
q_1(x) = q_2 (x) \text{ for a.e. } x\in [0,a].
\end{equation}
\end{theorem}

For reasons of brevity we stated Theorem~\ref{t1.1} only in 
the simplest possible case. 
Extensions to finite intervals $[0,R]$ instead of the 
half-line 
$[0,\infty)$, a 
discussion of boundary conditions other than Dirichlet at 
$x=0_+$, and the case of 
matrix-valued Schr\"odinger operators --- a new 
result --- will be
provided in the  main body of this paper. 

Theorem~\ref{t1.1} should be viewed as a local (and 
hence stronger) 
version of the 
following celebrated Borg-Marchenko uniqueness theorem, 
published 
by Marchenko 
\cite{Ma50} in 1950. Marchenko's extensive treatise 
on spectral 
theory of 
one-dimensional Schr\"odinger operators \cite{Ma52}, 
repeating 
the proof of his 
uniqueness theorem, then appeared in 1952, which 
also marked the 
appearance of 
Borg's proof of the uniqueness theorem \cite{Bo52} 
(apparently, 
based on his lecture 
at the 11th Scandinavian Congress of Mathematicians held at 
Trondheim, Norway in 
1949).

\begin{theorem} \lb{t1.2} \mbox{\rm (\cite{Bo52,Ma50,Ma52})} 
Suppose 
\begin{equation} \lb{1.3}
m_1(z) = m_2(z), \quad z\in\bbC\backslash\bbR,
\end{equation}
then
\begin{equation} \lb{1.4}
q_1(x) = q_2(x) \text{ for a.e. }  x\in [0,\infty).
\end{equation}
\end{theorem}

Again, we emphasize that Borg and Marchenko also treat the 
general case of non-Dirichlet 
boundary conditions at $x=0_+$, whose discussion we defer to 
Section~\ref{s2}. Moreover, 
Marchenko simultaneously discussed the half-line and finite 
interval case, also to be 
deferred to Section~\ref{s2}.

As pointed out by Levitan \cite{Le87} in the Notes to Chapter~2, 
Borg and
Marchenko  were actually preceded by Tikhonov \cite{Ti49} in 
1949, who
proved a special case of  Theorem~\ref{t1.2} in connection 
with the
string equation (and hence under certain  additional 
hypotheses on
$q_j$). Since Weyl-Titchmarsh functions $m(z)$ are uniquely  
related to
the spectral measure $d\rho$ of a self-adjoint (Dirichlet) 
Schr\"odinger 
operator $H=-\f{d^2}{dx^2} + q$ in $L^2 ([0,\infty))$ by the 
standard
Herglotz  representation
\begin{equation} \lb{1.5}
m(z) = \Real(m(i)) + \int_\bbR d\rho(\lambda) 
[(\lambda -z)^{-1} - 
\lambda(1+\lambda^2)^{-1}], \quad z\in\bbC\backslash\bbR,
\end{equation}
Theorem~\ref{t1.2} is equivalent to the following statement: 
Denote by $d\rho_j$ the 
spectral measures of $H_j$, $j=1,2$. Then
\begin{equation} \lb{1.6}
d\rho_1 = d\rho_2 \text{ implies } q_1 = q_2 \text{ a.e.~on } 
[0,\infty).
\end{equation}
In fact, Marchenko's proof takes the spectral measures 
$d\rho_j$ as the point of 
departure while Borg focuses on the Weyl-Titchmarsh functions 
$m_j$.

To the best of our knowledge, the only alternative approaches 
to Theorem~\ref{t1.2} 
are based on the Gelfand-Levitan solution of the inverse 
spectral problem published 
in 1951 (see also Levitan and Gasymov \cite{LG64}) and 
alternative 
variants due to 
M.~Krein \cite{Kr51}, \cite{Kr53}. In particular, it took over 
45 years to improve 
on Theorem~\ref{t1.2} and derive its local counterpart, 
Theorem~\ref{t1.1}. While the 
original proof of Theorem~\ref{t1.1} in \cite{Si98} relied on 
the full power of a new 
formalism in inverse spectral theory, relating $m(z)$ to finite 
Laplace transforms of 
the type
\begin{equation} \lb{1.7}
m(z) = iz^{1/2} - \int_0^a d\alpha\, A(\alpha) \, 
e^{2\alpha iz^{1/2}} + 
\widetilde O(e^{2\alpha iz^{1/2}})
\end{equation}
as $|z|\to\infty$ with $\arg(z)\in (\veps, \pi-\veps)$ for 
some $0<\veps<\pi$ (with 
$f=\widetilde O(g)$ if $g\to 0$ and for all $\delta >0$, 
$(\f{f}{g}) |g|^\delta\to 0$), 
we will present a short and fairly elementary argument in 
Section~\ref{s2}. In fact, 
as a corollary to our new proof of Theorem~\ref{t1.1}, we also 
obtain an elementary 
proof of a strengthened version of Theorem~\ref{t1.2}.

We should also mention some work of Ramm \cite{Ra99}, 
\cite{Ra99a}, who
provided a  proof of Theorem~\ref{t1.2} under a very strong 
additional
assumption, namely,  that $q_1$ and $q_2$ are both of short 
range. While
his result is necessarily weaker than the original 
Borg-Marchenko result,
Theorem~\ref{t1.2}, his method  of proof has elements in 
common with parts
of our proof (namely, he uses 
\eqref{2.27} below with $a=\infty$ and obtains a Volterra 
integral
equation close  to our \eqref{2.33}).

Finally, we have in preparation \cite{GS00} still
another  alternate proof of the local Borg-Marchenko theorem.

Extensions to finite intervals and general 
(i.e., non-Dirichlet) 
boundary conditions 
complete Section~\ref{s2}. Matrix-valued extensions of 
Theorem~\ref{t1.1} are presented in Section~\ref{s3}.

\section{A New Proof of Theorem~1.1} \lb{s2}

Throughout this section, unless explicitly stated otherwise, 
potentials
$q$ are supposed to satisfy 
\begin{equation} \lb{2.1}
q\in L^1 ([0,R]) \text{ for all } R>0, 
\quad \text{$q$ real-valued}.
\end{equation}
Given $q$, we introduce the corresponding self-adjoint 
Schr\"odinger operator $H$ in 
$L^2 ([0,\infty))$ with a Dirichlet boundary condition at 
$x=0_+$, by
\begin{align}
&H=-\f{d^2}{dx^2} + q, \no \\
&\dom(H) = \{g\in L^2 ([0,\infty)) \mid g,g' 
\in \AC([0,R])\text{ for all }R>0; \lb{2.2} \\
&\hspace*{16mm} g(0_+)=0, \,\sabc\text{ at } \infty; \, 
(-g''+qg) \in L^2
([0,\infty))\}. \no
\end{align}
Here ``s.-a.\,b.c." denotes a self-adjoint boundary condition 
at $\infty$ (which 
becomes relevant only if $q$ is in the limit circle case at 
$\infty$, but should 
be discarded otherwise, i.e., in the limit point case, where 
such a boundary 
condition is automatically satisfied). For example, an explicit 
form of such a 
boundary condition is
\begin{equation} \lb{2.3}
\lim_{x\uparrow \infty} W(f(z_0), g) (x) =0,
\end{equation}
where $f(z_0, x)$ for some fixed $z_0\in\bbC\backslash\bbR$, 
satisfies
\begin{equation} \lb{2.4}
f(z_0, \dott)\in L^2 ([0\infty)), \quad -f'' (z_0, x) 
+ [q(x)-z_0] 
f(z_0, x) =0 
\end{equation} 
and $W(f,g)(x) = f(x)g'(x) - f'(x)g(x)$ denotes the Wronskian 
of $f$ and $g$. 
Since these possible boundary conditions hardly play a role in 
the analysis 
to follow, we will not dwell on them any further. (Pertinent 
details can be 
found in \cite{GS96} and the references therein.)

Next, let $\psi(z,x)$ be the unique (up to constant multiples) 
Weyl solution associated 
with $H$, that is,
\begin{equation} \lb{2.5}
\begin{align}
&\psi(z, \dott) \in L^2 ([0,\infty)), 
\quad z\in\bbC\backslash\bbR, \notag \\
&\psi(z,x)\text{ satisfies the $\sabc$ of $H$ at $\infty$ 
(if any)}, \\
& -\psi'' (z,x) + [q(x) -z]\psi(z,x)=0. \notag
\end{align}
\end{equation}
Then the Weyl-Titchmarsh function $m(z)$ associated with $H$ 
is defined by
\begin{equation} \lb{2.6}
m(z) = \psi' (z,0_+) / \psi(z,0_+), 
\quad z\in\bbC\backslash\bbR
\end{equation}
and for later purposes we also introduce the corresponding 
$x$-dependent
version,
$m(z,x)$, by
\begin{equation} \lb{2.7}
m(z,x) = \psi'(z,x) / \psi(z,x), \quad z\in\bbC\backslash\bbR, 
\quad x\geq 0.
\end{equation}

After these preliminaries we are now ready to state the main 
ingredients used in 
our new proof of Theorem~\ref{t1.1}.

\begin{theorem}\lb{t2.1} \mbox{\rm(\cite{At81,Ev72})} Let 
$\arg(z)\in(\veps, \pi -\veps)$ 
for some $0<\veps < \pi$. Then for any fixed $x\in [0,\infty)$,
\begin{equation} \lb{2.8}
m(z,x) \eqlim_{|z|\to\infty} iz^{1/2} + o(1).
\end{equation}
\end{theorem}

The following result shows that one can also get an estimate 
uniform in $x$ as long as 
$x$ varies in compact intervals.

\begin{theorem}\lb{t2.2} \mbox{\rm(\cite{GS98})} Let 
$\arg(z)\in (\veps,\pi-\veps)$ for 
some $0<\veps<\pi$, and suppose $\delta>0$, $a>0$. Then there 
exists a $C(\veps,\delta,a) 
>0$ such that for all $x\in [0,a]$,
\begin{equation} \lb{2.9}
|m(z,x) - iz^{1/2}| \leq C(\veps,\delta, a),
\end{equation}
where $C(\veps,\delta,a)$ depends on $\veps,\delta$, and 
$\sup_{0\leq x\leq a} 
(\int_x^{x+\delta} dy\, |q(y)|)$.
\end{theorem}

Theorems~\ref{t2.1} and \ref{t2.2} can be proved following 
arguments of Atkinson 
\cite{At81}, who studied the Riccati-type equation satisfied 
by $m(z,x)$,
\begin{equation} \lb{2.10}
m'(z,x) + m(z,x)^2 = q(x)-z \text{ for a.e.}~x\geq 0 
\text{ and all }
z\in\bbC\backslash\bbR.
\end{equation}
Next, let $q_j(x)$, $j=1,2$ be two potentials satisfying 
\eqref{2.1}, 
with $m_j(z)$ the associated (Dirichlet) $m$-functions. 
Combining the 
{\it a priori} 
bound \eqref{2.9} with the differential equation resulting 
from \eqref{2.10},
\begin{align} 
&[m_1(z,x) - m_2(z,x)]' \lb{2.11}\\ 
& =q_1 (x) - q_2(x) -  [m_1 (z,x) + m_2(z,x)][m_1 (z,x) - m_2
(z,x)], \no
\end{align}
permits one to prove the following converse of 
Theorem~\ref{t1.1}.

\begin{theorem}\lb{t2.3} \mbox{\rm(\cite{GS98})} Let 
$\arg(z)\in (\veps, \pi -\veps)$ 
for some $0<\veps<\pi$ and suppose $a>0$. If
\begin{equation} \lb{2.12}
q_1(x) = q_2(x)\text{ for a.e. } x\in [0,a],
\end{equation}
then
\begin{equation} \lb{2.13}
|m_1 (z) - m_2(z)| \eqlim_{|z|\to\infty} 
O(e^{-2\Ima (z^{1/2})a}).
\end{equation}
\end{theorem}

\begin{lemma}\lb{l2.4} In addition to the hypotheses of 
Theorem~\ref{t2.2}, {\rm(}resp., 
Theorem~\ref{t2.3}{\rm)}, suppose that $H$ {\rm(}resp., 
$H_j$, $j=1,2${\rm)} is bounded 
from below. Then \eqref{2.9} {\rm(}resp., \eqref{2.13}{\rm)} 
extends to all $\arg(z) 
\in (\veps,\pi]$.
\end{lemma}

\begin{proof} Since $H_x \geq H$, where $H_x$ denotes the 
Schr\"odinger operator 
$-\f{d^2}{dx^2}+q$ in $L^2 ([x,\infty))$ with a Dirichlet 
boundary condition at $x_+$ 
(and the same $\sabc$ at $\infty$ as $H$, if any), there is 
an $E_0\in\bbR$ such that 
for all $x\in [0,a]$, $m(z,x)$ is analytic in 
$\bbC\backslash [E_0, \infty)$. Using 
$\ol{m(z,x)} = m(\bar z,x)$, the estimate \eqref{2.9} holds on 
the boundary of a 
sector with vertex at $E_0 -1$, symmetry axis 
$(-\infty, E_0 -1]$, and some opening 
angle $0<\veps<\pi/2$. An application of the 
Phragm\'en-Lindel\"of principle 
(cf.~\cite[Part~III, Sect.~6.5]{PS72}) then extends 
\eqref{2.9} to all of the 
interior of that sector and hence in particular along the 
ray $z\downarrow -\infty$. 
Since \eqref{2.13} results from \eqref{2.9} upon integrating 
(cf.~\eqref{2.11}),
\begin{align} 
& m_1(z,x) - m_2(z,x)]' \no \\
&=-[m_1(z,x) + m_2(z,x)][m_1(z,x), - m_2(z,x)], 
\quad x\in [0,a] \lb{2.14}
\end{align}
from $x=0$ to $x=a$, the extension of \eqref{2.9} to $z$ 
with $\arg(z)\in (\veps,\pi]$ 
just proven, allows one to estimate
\begin{equation} \lb{2.15}
|m_1(z,x) + m_2(z,x)| \eqlim_{|z|\to\infty} 2iz^{1/2} + O(1), 
\quad \arg(z)\in (\veps, \pi],
\end{equation}
uniformly with respect to $x\in [0,a]$, and hence to extend 
\eqref{2.13} to $\arg(z)\in 
(\veps, \pi]$. 
\end{proof}

\smallskip
Next, we briefly recall a few well-known facts on compactly 
supported $q$. Hence we 
suppose temporarily that
\begin{equation} \lb{2.15a}
\sup(\supp(q)) = \alpha <\infty.
\end{equation}
In this case, the Jost solution $f(z,x)$ associated with 
$q(x)$ satisfies
\begin{align}
f(z,x) &= e^{iz^{1/2}x} - \int_x^\alpha dy \, 
\f{\sin(z^{1/2}(x-y))}{z^{1/2}}\, 
q(y)\, f(z,y) \lb{2.16} \\
&= e^{iz^{1/2}x} + \int_x^\alpha dy K(x,y)\, e^{iz^{1/2}y}, 
 \quad  \Ima(z^{1/2}) \geq 0, \ x\geq 0, \lb{2.17}
\end{align}
where $K(x,y)$ denotes the transformation kernel satisfying 
(cf.~\cite[Sect.~3.1]{Ma86})
\begin{align}
K(x,y)&= \f12 \int_{(x+y)/2}^\alpha dx' \, q(x') 
- \int_{(x+y)/2}^\alpha \int_0^{(y-x)/2} 
dx'' \, q(x'-x'')\times \notag \\
&\hspace*{3.9cm} \times K(x'-x'', x' + x''), \quad x\leq y, 
\lb{2.18}
\\ K(x,y)&= 0, \quad x>y, \lb{2.19} \\
|K(x,y) &\leq \f12 \int_{(x+y)/2}^\alpha dx' \, |q(x')| 
\exp \left(\int_x^\alpha 
dx'' \, x''|q(x'')| \right). \lb{2.20} 
\end{align}
Moreover, $f(z,x)$ is a multiple of the Weyl solution, 
implying
\begin{equation} \lb{2.21}
m(z,x) = f'(z,x) / f(z,x), \quad z\in\bbC\backslash\bbR, 
\ x\geq 0,
\end{equation}
and the Volterra integral equation \eqref{2.16} immediately 
yields
\begin{align}
|f(z,x)| &\leq Ce^{-\Ima(z^{1/2})x}, \quad \Ima(z^{1/2}) 
\geq 0, \ x\geq 0, \lb{2.22} \\
f(z,x) &\eqlim_{\substack{|z|\to\infty \\ 
\Ima (z^{1/2})\geq 0}} 
 e^{iz^{1/2}x}(1+O(|z|^{-1/2}), \quad x\geq 0, \lb{2.23} 
\end{align}

Our final ingredient concerns the following result on finite 
Laplace transforms.

\begin{lemma}\lb{l2.5} \mbox{\rm($=$ Lemma~A.2.1 in 
\cite{Si98})} 
Let $g\in L^1 ([0,a])$ and 
assume that $\int_0^a dy\, g(y)e^{-xy} 
\underset{x\uparrow\infty}{=}
O(e^{-xa})$.  Then $g(y) = 0$ for a.e.~$ y\in [0,a]$.
\end{lemma}

Given these facts, the proof of Theorem~\ref{t1.1} now becomes 
quite simple.

\medskip
\begin{proof}[Proof of Theorem 1.1] By Theorem~\ref{t2.3} we 
may assume, without loss 
of generality, that $q_1$ and $q_2$ are compactly supported 
such that
\begin{equation} \lb{2.24}
\supp (q_j) \subseteq [0,a], \quad j=1,2,
\end{equation}
and by Lemma~\ref{l2.4} we may suppose that \eqref{1.1} holds 
along the ray 
$z\downarrow -\infty$, that is,
\begin{equation} \lb{2.25}
|m_1(z) - m_2(z)| \eqlim_{z\downarrow -\infty} 
O(e^{-2|z|^{1/2}a}).
\end{equation}
Denoting by $m_j(z,x)$ and $f_j(z,x)$ the $m$-functions and 
Jost solutions associated 
with $q_j$, $j=1,2$, integrating the elementary identity
\begin{equation} \lb{2.26}
\f{d}{dx}\, W(f_1(z,x), f_2(z,x)) = -[q_1(x) - q_2(x)] 
f_1(z,x) f_2(z,x)
\end{equation}
from $x=0$ to $x=a$, taking into account \eqref{2.21}, yields
\begin{align} 
&\int_0^a dx\, [q_1(x) - q_2(x)] f_1(z,x)f_2(z,x) \no \\
&= f_1(z,x) f_2(z,x) [m_1(z,x) - m_2(z,x)]\bigg|_{x=0}^a. 
\lb{2.27}
\end{align}
By \eqref{2.8}, \eqref{2.22}, and \eqref{2.25}, the right-hand 
side of \eqref{2.27} is 
$O(e^{-2|z|^{1/2}a})$ as $z\downarrow -\infty$, that is,
\begin{equation} \lb{2.28}
\int_0^a dx\, [q_1(x) - q_2(x)] f_1(z,x) f_2(z,x) 
 \eqlim_{z\downarrow -\infty} O(e^{-2|z|^{1/2}a}).
\end{equation}
Denoting by $K_j (x,y)$ the transformation kernels associated 
with $q_j$, $j=1,2$, 
\eqref{2.17} implies
\begin{equation} \lb{2.29}
f_1(z,x) f_2(z,x) = e^{2iz^{1/2}x} + \int_x^a dy \, L(x,y) 
\, e^{2iz^{1/2}y},
\end{equation}
where
\begin{align}
L(x,y) &= 2[K_1(x, 2y-x) + K_2(x, 2y-x)] \notag \\
&\quad + 2 \int_x^{2y-x} dx'\, K_1(x, x') K_2(x, 2y-x'), 
\quad x\leq y,
\lb{2.30} \\ L(x,y)&= 0, \quad x > y \quad \text{or} 
\quad y>a. \lb{2.31}
\end{align}
Insertion of \eqref{2.29} into \eqref{2.28}, interchanging 
the order of integration in 
the double integral, then yields
\begin{align} 
&\int_0^a dx [q_1 (x) - q_2(x)] f_1(z,x) f_2(z,x) \no \\
&= \int_0^a dy\left\{ [q_1(y) - q_2(y)] + \int_0^y dx \, L(x,y) 
[q_1(x) - q_2(x)]\right\} \, e^{-2|z|^{1/2}y} \no \\
& \eqlim_{z\downarrow -\infty} O(e^{-2|z|^{1/2}a}). \lb{2.32}
\end{align}
An application of Lemma~\ref{l2.5} then yields
\begin{equation} \lb{2.33}
[q_1(y) - q_2(y)] + \int_0^y dx\, L(x,y) [q_1(x) - q_2(x)] = 0 
\quad\text{for a.e. } y\in [0.a].
\end{equation}
Since \eqref{2.33} is a homogeneous Volterra integral 
equation with a continuous 
integral kernel $L(x,y)$, one concludes $q_1 = q_2$ 
a.e.~on $[0,a]$. 
\end{proof}
In particular, one obtains the following strengthened version 
of the original 
Borg-Marchenko uniqueness result, Theorem~\ref{t1.2}.

\begin{corollary}\lb{c2.6} Let $0<\veps <\pi/2$ and suppose 
that for all $a>0$,
\begin{equation} \lb{2.34}
|m_1(z) - m_2(z)| \eqlim_{|z|\to\infty} 
O(e^{-2\Ima (z^{1/2})a})
\end{equation}
along the ray $\arg(z) = \pi -\veps$. Then
\begin{equation} \lb{2.35}
q_1(x) = q_2(x) \text{ for a.e. } x\in [0,\infty).
\end{equation}
\end{corollary}

\begin{remark} \lb{r2.7} The Borg-Marchenko uniqueness result, 
Theorem~\ref{t1.2} 
(but not our strengthened version, Corollary~\ref{c2.6}), 
under the
additional  condition of short-range potentials $q_j$ 
satisfying 
$q_j\in L^1 ([0,\infty); 
(1+x)\, dx)$, $j=1,2$, can also be proved using Property~C, 
a device
recently used  by Ramm \cite{Ra99,Ra99a} in a variety of 
uniqueness
results. In this case, \eqref{2.27}  for $z=\lambda>0$ becomes
\begin{align} 
&\int_0^\infty dx\, [q_1(x) - q_2(x)] f_1(\lambda,x) 
f_2(\lambda, z) \no
\\ 
&= -f_1(\lambda, 0) f_2(\lambda, 0) [m_1(\lambda +i0) -
m_2(\lambda+i0)] = 0, 
\quad \lambda >0 \lb{2.36}
\end{align}
since $m_1(z) = m_2(z)$, $z\in\bbC_+$ extends to 
$m_1(\lambda+i0) = m_2(\lambda+i0)$, 
$\lambda >0$ by continuity in the present short-range case. 
By definition, Property~C 
stands for completeness of the set 
$\{f_1(\lambda, x)f_2(\lambda, x)\}_{\lambda>0}$ 
in $L^1 ([0,\infty); (1+x)\, dx)$ (this extends to 
$L^1 ([0,\infty))$)
and hence
\eqref{2.36}  yields $q_1 = q_2$ a.e.~on $[0,\infty)$.
\end{remark}
In the remainder of this section, we consider a variety of 
generalizations of the 
result obtained.

\begin{remark} \lb{r2.8} The ray $\arg(z) = \pi -\veps$, 
$0<\veps < \pi/2$ chosen in 
Theorem~\ref{t1.1} and Corollary~\ref{2.6} is of no particular 
importance. A limit 
taken along any non-self-intersecting curve $\calC$ going to 
infinity in the sector 
$\arg(z)\in (\pi/2+\veps, \pi -\veps)$ will do as we can apply 
the Phragm\'en-Lindel\"of 
principle (\cite[Part~III, Sect.~6.5]{PS72}) to the region 
enclosed by $\calC$ and its 
complex conjugate $\bar\calC$ (needed in connection with 
Lemma~\ref{l2.4} in order to 
reduce the general case to the case of spectra bounded from 
below).
\end{remark}

\begin{remark} \lb{r2.9} For simplicity of exposition, we only 
discussed the 
Dirichlet boundary condition
\begin{equation} \lb{2.37}
g(0_+)=0
\end{equation}
in the definition of $H$ in \eqref{2.2}. Next we replace 
\eqref{2.37} by the general 
boundary condition
\begin{equation} \lb{2.38}
\sin(\alpha) g'(0_+) + \cos(\alpha) g(0_+) = 0, \quad 
\alpha\in [0,\pi)
\end{equation}
in \eqref{2.2}, denoting the resulting Schr\"odinger operator 
by $H_\alpha$, while 
keeping the boundary condition at infinity (if any) identical 
for all $\alpha\in 
[0,\pi)$. Denoting by $m_\alpha (z)$ the Weyl-Titchmarsh 
function associated with 
$H_\alpha$, the well-known relation (cf.~e.g., Appendix~A of 
\cite{GS96} for precise 
details on $H_\alpha$ and $m_\alpha(z)$)
\begin{equation} \lb{2.39}
m_\alpha(z) = \f{-\sin(\alpha) + \cos(\alpha) m(z)}{\cos(\alpha) 
+ \sin(\alpha) m(z)}\, , 
\quad \alpha\in [0,\pi), \ z\in \bbC\backslash\bbR
\end{equation}
reduces the case $\alpha\in (0,\pi)$ to the Dirichlet case 
$\alpha =0$.
In  particular, Theorem~\ref{t1.1} and Corollary~\ref{c2.6} 
remain valid
with $m_j (z)$  replaced by $m_{j,\alpha}(z)$, 
$\alpha\in [0,\pi)$.
Indeed, $|m_{1,\alpha}(z) -  m_{2,\alpha}(z)| 
\underset{|z|\to\infty}{=}
O(e^{-2\Ima (z^{1/2})a})$ along the  ray $\arg(z) = 
\pi -\veps$ is easily
seen to imply, for all sufficiently small $\delta >0$,
\begin{equation} \lb{2.40}
\begin{aligned}
|m_{1,0}(z) - m_{2,0}(z)| &\eqlim_{|z|\to\infty} O(|z|\, 
e^{-2\Ima(z^{1/2})a}) \\
&\eqlim_{|z|\to\infty} O(e^{-2\Ima (z^{1/2})(a-\delta)})
\end{aligned}
\end{equation}
along the ray $\arg(z) = \pi - \veps$. Hence one infers 
from Theorem~\ref{t1.1} 
that for all $0<\delta <a$, $q_1 = q_2$ a.e.~on $[0, 
a-\delta]$. Since
$\delta > 0$  can be chosen arbitrarily small, one 
concludes $q_1 = q_2$
a.e.~on $[0,a]$. In fact,  more is true. Since
$m_\alpha(z)\underset{|z|\to\infty}{\to}\cot(\alpha)$ 
along the ray, one
concludes that 
$|m_{1,\alpha_1}  - m_{2,\alpha_2}(z)| 
\underset{|z|\to\infty}{=}
O(e^{-2\Ima (z^{1/2})a})$ along a ray implies $\alpha_1 
= \alpha_2$
and $q_1 = q_2$ a.e.~on $[0, a]$.
\end{remark}

\begin{remark} \lb{r2.10} If one is interested in a finite 
interval
$[0,b]$ instead  of the half-line $[0,\infty)$ in 
Theorem~\ref{t1.1}, with
$0<a<b$, one introduces  a self-adjoint boundary condition 
at $x=b_-$ of
the type
\begin{equation} \lb{2.41}
\sin(\beta) g'(b_-) + \cos(\beta) g(b_-) = 0, 
\quad \beta\in [0,\pi).
\end{equation}
The analog of the Weyl solution $\psi (z,x; b,\beta)$ for 
the corresponding 
Schr\"odinger operator $H(b,\beta)$ in $L^2 ([0,b])$ 
defined by
\begin{align}
 &H (b,\beta) = -\f{d^2}{dx^2} + q, \quad \beta\in [0,\pi), 
\lb{2.42} \\
& \dom(H(b,\beta)) = \{g\in L^2([0,b]) 
\mid g,g'\in \AC([0,b]); 
\ g(0_+)=0, \no \\
&\hspace*{9mm} \sin(\beta) g'(b_-) + \cos(\beta) g(b_-) =0; 
\ (-g'' + qg)\in L^2 ([0,b])\} \no
\end{align}
is then defined by 
\begin{align} 
&\sin(\beta)\psi'(z, b_-; b,\beta) + \cos(\beta) 
\psi(z, b_-; b,\beta) =
0, \quad  z\in\bbC\backslash\bbR, \no \\
&-\psi''(z,x; b,\beta) + [q(x) - z] \psi(z,x; b,\beta)=0. 
\lb{2.43}
\end{align}
Moreover, the analog of \eqref{2.17} is then of the type
\begin{multline} \lb{2.44}
\psi(z,x; b,\beta) = \psi^{(0)} (z,x;b,\beta) 
 + \int_x^b dy \, 
K(x,y; b,\beta)\psi^{(0)}(z,x; b,\beta), \\
\quad z\in\bbC\backslash\bbR, \ x\in [0,b],
\end{multline}
where
\begin{align}
&\psi^{(0)} (z,x;b,\beta) = \f{e^{iz^{1/2}x} 
+ \zeta(\beta,z)\,
e^{iz^{1/2}(2b-x)}} {1 + \zeta(\beta, z)\, 
e^{2iz^{1/2}b}},\, \lb{2.45} \\
&\zeta(\beta, z) = \f{-iz^{1/2} - \cot(\beta)}{-iz^{1/2} 
+ \cot(\beta)}\,
, \quad \Ima(z^{1/2})\geq 0 \no
\end{align}
is the corresponding Weyl solution in the case $q(x)=0$, 
$x\in [0,b]$, and $K(x,y;b,\beta)$ 
is a transformation kernel analogous to $K(x,t;h)$ discussed 
in Sect.~1.3 of \cite{Ma86}. 
Theorem~\ref{t1.1} then extends to triples 
$(q_j, b_j, \beta_j)$, $j=1,2$ with $a<\min 
(b_1, b_2)$, replacing $f_j (z,x)$ in \eqref{2.28} by 
$\psi(z,x; b_j, \beta_j)$, 
$j=1,2$. More precisely, if $m(z; b_j, \beta_j)$ denote the 
$m$-functions for $H(b_j, 
\beta_j)$, $j=1,2$ with $a<\min (b_1, b_2)$ and
\begin{equation} \lb{2.46}
|m(z; b_1, \beta_1) - m(z; b_2, \beta_2)| 
\eqlim_{|z|\to\infty} O(e^{-2\Ima(z^{1/2})a}) 
\end{equation}
along the ray $\arg(z) = \pi -\veps$, then $q_1 = q_2$ 
a.e.~on $[0,a]$. 
In
fact, it was precisely this version of Theorem~\ref{t1.1} 
which was
originally proven  by one of us \cite{Si98} in 1998. 

One can also derive
additional results in the case 
$a=b_1 = b_2$ (cf.~Theorem~1.3 in \cite{Si98}). Indeed, 
$\zeta(\beta, z) \underset{|z|\to\infty}{=} 1 + 
O(|z|^{-1/2})$, so by \eqref{2.44} and \eqref{2.45},
\begin{align}
\psi (z,0; a,\beta) &\underset{|z|\to\infty}{=} 
1 + O(|z|^{-1/2}),
\lb{new2.47} \\
\psi (z,a; a,\beta) &\underset{|z|\to\infty}{=} 
2e^{iz^{1/2}a}
(1+O(|z|^{-1/2})), \lb{new2.48}
\end{align}
which are analogous to \eqref{2.23}.
Thus, if 
$q_1 = q_2$ on $[0,a]$ but $\beta_1 \neq \beta_2$, we have 
that
\begin{align}
&[m_1 (z; a, \beta_1) - m_2 (z; a,\beta_2)] \lb{new2.50}  \\
&= \psi_1 (z,0; a,\beta_1)^{-1} \psi_2 (z, 0; a, \beta_2)^{-1} 
W(\psi_2 (z,0;a,\beta_2), \psi_1 (z,0; a,\beta_1)) \no \\
&= \psi_1 (z;0, a, \beta_1)^{-1} \psi_2 (z,0; a,\beta_2)^{-1} 
W(\psi_2 (z,a; a,\beta_2), \psi_1 (z,a; a,\beta_2)), \no
\end{align}
by the constancy of the Wronskian when $q_1 = q_2$. But 
$\psi'(a)/\psi(a)$ equals $\cot(\beta)$ by \eqref{2.43} and 
hence 
\begin{multline} \lb{new2.51}
m_1 (z; a,\beta_1) - m_2 (z; a, \beta_2) = 
\psi_1 (z, 0; a, \beta_1)^{-1} 
\psi_2 (z, 0; a,\beta_2)^{-1}\times \\
\qquad \times\psi_1 (z,a; a,\beta_1) \psi_2 (z,a; a,\beta_2) 
[\cot(\beta_1) - \cot(\beta_2)].
\end{multline}
Using \eqref{new2.47}, \eqref{new2.48}, this implies that
\begin{align}
&m_1 (z; a,\beta_1) - m_2 (z; a,\beta_2) \no \\
&= 4e^{2iz^{1/2}a} 
[\cot(\beta_1) - \cot(\beta_2)][1+O(|z|^{-1/2})], \lb{2.51a}
\end{align}
which is Theorem~1.3 in \cite{Si98}.
\end{remark}


While we have separately described a few extensions in 
Remarks~\ref{r2.8}--\ref{r2.10}, 
it is clear that they can all be combined at once.

We also mention the analog of Theorem~\ref{t1.1} for 
Schr\"odinger
operators on the real line. Assuming 
\begin{equation} \lb{2.46a}
q\in L^1_{\loc} (\bbR), \quad \text{$q$ real-valued,}
\end{equation}
one introduces the corresponding self-adjoint Schr\"odinger 
operator $H$
in
$L^2 (\bbR)$ by
\begin{align}
&H=-\f{d^2}{dx^2} + q, \lb{2.46b} \\
&\dom(H) = \{g\in L^2 (\bbR) \mid g,g' 
\in \AC_{\loc}(\bbR); \,\text{s.}\,\sabc\text{ at } 
\pm \infty; \no \\
& \hspace*{7.3cm}   
(-g''+qg) \in L^2 (\bbR)\}. \no
\end{align}
Here ``s.\,s.-a.\,b.c." denotes separated self-adjoint 
boundary conditions
at $+\infty$ and/or $-\infty$ (if any). 

The $2\times 2$ matrix-valued $m$-function $\calM (z)$ 
associated with $H$
in $L^2(\bbR)$ is then defined by 
\begin{align}
\calM (z) &=(m_-(z)-m_+(z))^{-1} \times \lb{2.46c} \\ 
& \quad \times \begin{pmatrix}
1& (m_-(z)+m_+(z))/2 \\ (m_-(z)+m_+(z))/2
& m_-(z)m_+(z) \end{pmatrix}, 
\quad z\in\bbC\backslash\bbR, \no
\end{align}
where $m_\pm (z)$ denote the half-line $m$-functions 
associated with $H$ 
restricted to $[0,\pm\infty)$ and a Dirichlet boundary 
condition at
$x=0$. 

Next, let $q_j(x)$, $j=1,2$ be two potentials satisfying 
\eqref{2.46a}
and $H_j$ the corresponding Schr\"odinger operators 
\eqref{2.46b} in
$L^2(\bbR)$, with $\calM_j (z)$, $j=1,2$ the associated 
$2\times 2$ 
matrix-valued 
$m$-functions. Then the analog of Theorem~\ref{t1.1} reads 
as follows.

\begin{theorem}\lb{t2.11} Let $a>0$, $0<\veps<\pi/2$ and 
suppose that
\begin{equation} \lb{2.46d}
\|\calM_1 (z) - \calM_2(z)\|_{\bbC^{2\times 2}} 
\eqlim_{|z|\to\infty}
O(e^{-2\Ima (z^{1/2})a}) 
\end{equation} 
along the ray $\arg(z) = \pi-\veps$. Then
\begin{equation} \lb{2.46e}
q_1(x) = q_2 (x) \text{ for a.e. } x\in [-a,a].
\end{equation}
\end{theorem}

\begin{proof}
We denote by $m_{j,\pm}(z)$ the half-line (Dirichlet) 
$m$-functions
associated with $H_j$ on $[0,\pm\infty)$, $j=1,2$. Then a 
straightforward combination of \eqref{2.8} and 
\eqref{2.46d} yields
\begin{equation}
|m_{1,\pm}(z)-m_{2,\pm}(z)|\underset{|z|\to\infty}{=}
O(|z|e^{-2\Ima (z^{1/2})a}) \lb{2.46f}
\end{equation}
and hence \eqref{2.46e}, applying Theorem~\ref{t1.1} 
separately 
to the two half-lines $[0,\infty)$ and $(-\infty,0]$ 
(and using the
argument  following \eqref{2.40}.
\end{proof}

Finally, the reader might be interested in the analog of 
Theorem~\ref{t1.1} in the case of 
second-order difference operators, that is, Jacobi operators. 
Let $A$ be a bounded 
self-adjoint Jacobi operator in $\ell^2 (\bbN_0)$ ($\bbN_0 
= \bbN\cup\{0\}$) of the 
type
\begin{equation} \lb{2.47}
A = \left(
\begin{matrix}
b_0 & a_0 & 0 & 0 & \dots & \dots\\
a_0 & b_1 & a_1 & 0 & \dots & \dots \\
0 & a_1 & b_2 & a_2 & \dots & \dots\\
0 & 0 & a_2 & \ddots & \ddots \\
\vdots & \vdots & \vdots & \ddots & \ddots & \ddots \\
\vdots & \vdots & \vdots & {} & \ddots & \ddots
\end{matrix} \right), \quad a_k >0, \ b_k\in\bbR, \ k\in\bbN_0.
\end{equation}
The corresponding $m$-function of $A$ is then defined by
\begin{equation} \lb{2.48}
m(z) = (\delta_0, (A-z)^{-1}\delta_0) 
= \int_\bbR d\rho(\lambda) (\lambda -z)^{-1}, 
\quad z\in\bbC\backslash\bbR,
\end{equation}
where $\delta_0 = (1,0,0,\dots)$. The analog of 
Theorem~\ref{t1.1} in the discrete case 
then reads as follows. Denote by $m_j(z)$ the $m$-functions 
for two self-adjoint Jacobi 
operators $A_j$, $j=1,2$, denoting the matrix elements of 
$A_j$ by $a_{j,k}$, $b_{j,k}$, 
$j=1,2$, $k\in \bbN_0$. Then
\begin{equation} \lb{2.49}
|m_1(z) - m_2(z)| \eqlim_{|z|\to\infty} O(|z|^{-N}),
\end{equation}
for some $N\in\bbN$, $N\geq 3$, if and only if
\begin{equation} \lb{2.50}
a_{1,k} = a_{2,k}, \quad b_{1,k} = b_{2,k}, 
\quad 0\leq k \leq \f{N-4}{2} 
\quad\text{if $N$ is even} \quad (N\geq 4)
\end{equation}
and
\begin{equation} \lb{2.51}
\begin{aligned}
a_{1,k} &= a_{2,k}, \quad 0\leq k \leq \f{N-5}{2}\, , \\ 
b_{1,k} &= b_{2,k}, \quad 0\leq k \leq \f{N-3}{2} 
\quad\text{if $N$ is odd}.
\end{aligned}
\end{equation}
The proof is clear from \eqref{2.48} and the well-known 
formulas 
(cf.~\cite[Sect.~VII.1]{Be68}).
\begin{equation} \lb{2.52}
a_k = \int_\bbR d\rho(\lambda) \, \lambda P_k(\lambda) 
P_{k+1}(\lambda), 
\quad b_k = \int_\bbR d\rho(\lambda)\, \lambda P_k(\lambda)^2, 
\quad 
k\in\bbN_0,
\end{equation}
where $\{P_k(\lambda)\}_{k\in\bbN_0}$ is an orthonormal 
system of polynomials with 
respect to the spectral measure $d\rho$, with $P_k(z)$ of 
degree $k$ in $z$, $P_0(z) =1$.

\section{Matrix-Valued Schr\"odinger Operators} \lb{s3}

In our final section we extend Theorem~\ref{t1.1} to 
matrix-valued potentials 
(cf., \cite[Ch.~III]{CL90}, \cite{KS88}, \cite{La84} and 
the references therein).

Let $m\in\bbN$ and denote by $I_m$ the identity matrix 
in $\bbC^m$. Assuming
\begin{equation} \lb{3.1}
Q=Q^* \in L^1 ([0,R])^{m\times m} \text{ for all } R>0,
\end{equation}
we introduce the corresponding matrix-valued self-adjoint 
Schr\"odinger operator 
$H$ in $L^2 ([0,\infty))^m$ with a Dirichlet boundary 
condition at $x=0_+$, by 
\begin{equation} 
\begin{align}
&H=-\f{d^2}{dx^2}\,I_m + Q, \lb{3.2} \\
&\dom(H) = \{g\in L^2 ([0,\infty))^m \mid g,g' 
\in \AC([0,R])^m\text{ for all }R>0; \no \\
&\hspace*{17mm}  g(0_+)=0, \,\sabc\text{ at } 
\infty; \, (-g''+Qg) 
\in L^2 ([0,\infty))^m\}. \no
\end{align}
\end{equation}
Here ``$\sabc$ at $\infty$" again denotes a self-adjoint 
boundary condition at $\infty$ 
(if $Q$ is not in the limit point case at $\infty$). For 
more details about the limit 
point/limit circle and all the intermediate cases, see 
\cite{CG99,HS81,HS83,HS84,KR74,
Kr89a,Kr89b,Or76,Ro69} and the references therein.

Next, let $\Psi(z,x)$ be the unique (up to right multiplication 
of non-singular constant 
$m\times m$ matrices) $m\times m$ matrix-valued Weyl solution 
associated with $H$, 
satisfying
\begin{align}
&\Psi(z,\dott) \in L^2 ([0,\infty))^{m\times m}, \quad
z\in\bbC\backslash\bbR, \lb{3.3} \\
&\Psi(z,x)\text{ satisfies the $\sabc$ of $H$ at 
$\infty$ (if any)},
\lb{3.4} \\ 
&-\Psi''(z,x) + [Q(x) - zI_m]\Psi(z,x) = 0. \lb{3.5}
\end{align}

The $m\times m$ matrix-valued Weyl-Titchmarsh function 
$M(z)$ associated with $H$ is then 
defined by
\begin{equation} \lb{3.6}
M(z) = \Psi' (z, 0_+) \Psi(z, 0_+)^{-1}, 
\quad z\in\bbC\backslash\bbR 
\end{equation}
and similarly, we introduce its $x$-dependent version, 
$M(z,x)$, by
\begin{equation} \lb{3.7}
M(z,x) = \Psi' (z,x) \Psi(z,x)^{-1}, 
\quad z\in\bbC\backslash\bbR, \, x\geq 0.
\end{equation}
The matrix Riccati equation satisfied by $M(z,x)$, the 
analog of
\eqref{2.10}, then  reads
\begin{equation} \lb{3.8}
M'(z,x) + M(z,x)^2 = Q(x) - zI_m \text{ for a.e.}~x\geq 0 
\text{ and all } z\in\bbC\backslash\bbR.
\end{equation}
Next, let $Q_j(x)$, $j=1,2$ be two self-adjoint matrix-valued 
potentials satisfying 
\eqref{3.1}, and $M_j(z)$, $M_j(z,x)$ the Weyl-Titchmarsh 
matrices associated with the 
corresponding (Dirichlet) Schr\"odinger operators. Then 
the analog of \eqref{2.11} is 
of the form
\begin{align}
&[M_1(z,x) - M_2(z,x)]'  \lb{3.9} \\ 
&=Q_1 (x) - Q_2(x) - \tfrac{1}{2} [M_1 (z,x) + M_2(z,x)]
[M_1 (z,x) - M_2 (z,x)] \no \\
& \quad  - \tfrac{1}{2} [M_1 (z,x) 
- M_2(z,x)][M_1 (z,x) + M_2
(z,x)]. \no
\end{align}
Combining \eqref{3.9} with the elementary fact that any 
$m\times m$ matrix-valued 
solution $U(x)$ of
\begin{equation} \lb{3.10}
U'(x)=B(x) U(x) + U(x)B(x)
\end{equation} 
is of the form
\begin{equation} \lb{3.11}
U(x) = V(x) CW(x),
\end{equation}
where $C$ is a constant $m\times m$ matrix and $V(x)$, 
respectively,
$W(x)$, is a  fundamental system of solutions of $R'(x) 
= B(x)R(x)$,
respectively, $S'(x) = S(x)B(x)$,  one can prove the 
analogs of
Theorems~\ref{t2.1}--\ref{t2.3} in the present matrix  
context. More
precisely, the matrix analogs of Theorems~\ref{t2.1} and 
\ref{t2.2} 
follow from Theorem~4.8 in \cite{CG99}. The corresponding 
analog of
Theorem~\ref{t2.3}  follows from Theorem~4.5 and 
Remark~4.7 in
\cite{CG99}. Moreover, in the case that $H$  is bounded 
from below,
Lemma~\ref{l2.4} generalizes to the matrix-valued context 
and  hence
permits one to take the limit $z\downarrow -\infty$ in 
the matrix
analog of 
\eqref{2.13}. While the scalar case treated in detail in 
\cite{GS98} is based on 
Riccati-type identities such as \eqref{2.11} and an 
{\it a priori} bound
of the type 
\eqref{2.9} inspired by Atkinson's 1981 paper \cite{At81}, 
the matrix-valued case 
discussed in depth in \cite{CG99} is based on corresponding 
Riccati-type identities 
such as \eqref{3.9} and an {\it a priori} bound of the type
\begin{equation} \lb{3.12}
M(z,x) = iz^{1/2} I_m + o(|z|^{1/2})
\end{equation}
first obtained by Atkinson in an unpublished manuscript 
\cite{At??}.

In the special case of short-range matrix-valued potentials 
$Q(x)$, $m\times m$ matrix 
analogs of the Jost solution $F(z,x)$ as well as the 
transformation kernel $K(x,y)$ 
associated with $H$ as in \eqref{2.16}--\eqref{2.20} 
(replacing $|\dott|$ by an 
appropriate matrix norm $\|\dott\|_{\bbC^{m\times m}}$, 
have been
discussed in great detail in the  classical 1963 monograph 
by Agranovich
and Marchenko
\cite[Ch.~I]{AM63}. Moreover, 
\eqref{2.21}--\eqref{2.23} trivially extend to the matrix case.

Given these preliminaries, the analog of Theorem~\ref{t1.1} 
and Corollary~\ref{c2.6} 
reads as follows in the matrix-valued context.

\begin{theorem}\lb{t3.1} Let $a>0$, $0<\veps<\pi/2$ and 
suppose
\begin{equation} \lb{3.13}
\|M_1(z) - M_2(z)\|_{\bbC^{m\times m}} 
\eqlim_{|z|\to\infty} O(e^{-2\Ima
(z^{1/2})a})
\end{equation}
along the ray $\arg(z) = \pi - \veps$. Then
\begin{equation} \lb{3.14}
Q_1(x) = Q_2(x) \text{ for a.e. } x\in [0,a].
\end{equation}
In particular, if \eqref{3.13} holds for all $a>0$, 
then $Q_1 = Q_2$ a.e.~on $[0, \infty)$.
\end{theorem}

\begin{proof}[Sketch of Proof] As in the scalar case, we 
may assume without loss of 
generality that
\begin{equation} \lb{3.15}
\supp (Q_j) \subseteq [0,a], \quad j=1,2.
\end{equation}
The fundamental identity \eqref{2.26}, in the present 
non-commutative case, needs to 
be replaced by
\begin{equation} \lb{3.16}
\f{d}{dx}\, W(F_1 (\bar z, x)^*, F_2(z,x)) = 
-F_1(\bar z,x)^* [Q_1(x) - Q_2(x)] 
F_2 (z,x),
\end{equation}
where $F_j(z,x)$ denote the $m\times m$ matrix-valued Jost 
solutions associated with 
$Q_j$, $j=1,2$, and $W(F,G)(x) = F(x) G'(x) - F'(x)G(x)$ 
the matrix-valued 
Wronskian of $m\times m$ matrices $F$ and $G$. Identity 
\eqref{2.27} then becomes
\begin{align} 
&\int_0^a dx\, F_1(\bar z,x)^* [Q_1(x) - Q_2(x)] 
F_2 (z,x) \no \\
&= F_1(\bar z, x)^* [M_1(z,x) - M_2(z,x)] 
F_2(z,x)\bigg|_{x=0}^a\, , 
\lb{3.17}
\end{align}
utilizing the fact
\begin{equation} \lb{3.18}
M_1(\bar z,x)^* = M_1(z,x).
\end{equation}
$F_j$ obeys a transformation kernel representation
\begin{align} 
&F_j (z,x) = e^{iz^{1/2}x} I_m + \int_x^a dy\, K_j (x,y)\, 
e^{iz^{1/2}y} I_m \, , \lb{3.19} \\
& \hspace*{3cm} \Ima(z^{1/2}) \geq 0, \, x\geq 0, \, 
j=1,2. \no
\end{align}
 From this, \eqref{3.12}, and the hypothesis of \eqref{3.13}, 
one concludes by \eqref{3.17} that
\begin{equation} \lb{3.20}
\int_0^a dx\, F_1 (\bar z, x) [Q_1(x) - Q_2(x)] F_2 (z,x) 
\underset{z\downarrow -\infty}{=} O(e^{-2\Ima(z^{1/2})a}).
\end{equation}

Now let $R_A$ be right multiplication by $A$ on $n\times n$ 
matrices and $L_B$ be 
left multiplication by $B$. Then
\begin{align} 
&\text{LHS of \eqref{3.20}} \lb{3.21} \\
&= \int_0^a dx\, \left\{Q_1(y) - Q_2(y) + \int_0^y dx\, 
\calL(x,y) [Q_1(x) - Q_2(x)]\right\} e^{-2z^{1/2}y}, \no
\end{align}
where $\calL$ is an operator on $n\times n$ matrices which 
is a sum of a left 
multiplication (by $2K_j(x, 2y-x)$), a right multiplication 
(by $2K_2 (x, 2y-x)$), and 
a convolution of a left and right multiplication.

It follows by Lemma~\ref{l2.5}, \eqref{3.20}, and 
\eqref{3.21} that
\begin{equation} \lb{3.22}
Q_1(y) - Q_2(y) + \int_0^y dx \, \calL(x,y) [Q_1 (x) 
- Q_2(x)] = 0.
\end{equation}

This is a Volterra equation and the same argument based on 
\[
\int_0^y dx_1 \int_0^{x_1} dx_2 \cdots \int_0^{x_{n-1}} 
dx_n = \f{y^n}{n!}
\]
that a Volterra operator has zero spectral radius applies to
operator-valued  Volterra equations. Thus, \eqref{3.22} 
mplies 
$Q_1(y) - Q_2(y) = 0$  for a.e.~$y\in [0,a]$.
\end{proof}

Extensions of Theorem~\ref{t3.1} in the spirit of
Remarks~\ref{r2.8}--\ref{r2.10} and Theorem~\ref{t2.11} 
can be made, but
we omit the corresponding details at this point.

\vspace*{2mm}
\noindent {\bf Acknowledgments.} F.~G.~thanks 
T.~Tombrello  for the hospitality of Caltech where this work 
was done.


\end{document}